\documentclass[a4paper,twoside,11pt]{article}

\usepackage{amsmath} \usepackage{amssymb} \usepackage{amsopn}
\usepackage{amsthm} \usepackage{amsfonts} \usepackage{mathrsfs}

\usepackage[dvips]{graphicx} \usepackage[usenames]{color}

\usepackage[all]{xy}

\usepackage[latin1]{inputenc}  

\usepackage[pdfpagelabels, dvips]{hyperref}

\usepackage{a4wide}
\newcommand{\Ueberschrift}{On the period-index problem in light of the section conjecture} 

\hypersetup{
  pdftitle={},
  pdfauthor={},
  pdfsubject={},
  pdfkeywords={},
  colorlinks=false,
  breaklinks=true,
  bookmarksopen=true,
  bookmarksnumbered=true,
  pdfpagemode=UseOutlines,
  plainpages=false}

\DeclareMathOperator{\rC}{C}

\DeclareMathOperator{\rH}{H}

\newcommand{\bC}{{\mathbb C}}

\newcommand{\bG}{{\mathbb G}}
\newcommand{\bH}{{\mathbb H}}

\newcommand{\bP}{{\mathbb P}}
\newcommand{\bQ}{{\mathbb Q}}
\newcommand{\bR}{{\mathbb R}}

\newcommand{\bZ}{{\mathbb Z}}

\newcommand{\cA}{{\mathscr A}}

\newcommand{\cE}{{\mathscr E}}

\newcommand{\cL}{{\mathscr L}}
\newcommand{\cM}{{\mathscr M}}

\newcommand{\dO}{{\mathcal O}}

\newcommand{\fo}{{\mathfrak o}}
\newcommand{\fp}{{\mathfrak p}}

\newcommand{\surj}{\twoheadrightarrow} 
\newcommand{\inj}{\hookrightarrow}

\DeclareMathOperator{\id}{id}

\newcommand{\ev}{{\rm ev}}

\DeclareMathOperator{\Aut}{Aut}
\DeclareMathOperator{\End}{End}

\DeclareMathOperator{\Mat}{M}

\DeclareMathOperator{\Char}{char} 

\DeclareMathOperator{\Spec}{Spec}

\DeclareMathOperator{\Div}{Div}
\DeclareMathOperator{\divisor}{div}
\DeclareMathOperator{\Pic}{Pic}
\DeclareMathOperator{\uPic}{\underline{Pic}}

\DeclareMathOperator{\Sing}{Sing}

\newcommand{\redu}{{\rm red}}
\newcommand{\norm}{{\rm norm}}

\DeclareMathOperator{\CH}{CH}

\newcommand{\OO}{\dO}

\DeclareMathOperator{\mg}{\cM_g}

\newcommand{\Gm}{\bG_m}

\DeclareMathOperator{\ix}{index}
\DeclareMathOperator{\pe}{period}

\DeclareMathOperator{\Br}{Br}
\newcommand{\inv}{{\rm inv}}
\DeclareMathOperator{\Frob}{Frob}

\DeclareMathOperator{\BS}{BS}

\DeclareMathOperator{\R}{R}

\DeclareMathOperator{\cAut}{\cA\mathit{ut}}

\DeclareMathOperator{\Tot}{Tot}

\DeclareMathOperator{\cores}{cor}

\DeclareMathOperator{\Gal}{Gal}

\newcommand{\ph}{\varphi}

\newcommand{\alg}{{\rm alg}}

\newcommand{\bruch}[2]{\genfrac{}{}{0.5pt}{}{#1}{#2}}

\newcommand{\ov}[1]{\mbox{${\overline{#1}}$}} 

\newtheorem{thm}{Theorem}

\newtheorem{prop}[thm]{Proposition}
\newtheorem{lem}[thm]{Lemma}

\newtheorem{cor}[thm]{Corollary}
\newtheorem{conj}[thm]{Conjecture}

\theoremstyle{definition}
\newtheorem{defi}[thm]{Definition}

\theoremstyle{remark}
\newtheorem{rmk}[thm]{Remark}

\newtheorem{ex}[thm]{Example}

\newenvironment{pro}[1][Proof]{{\it{#1:}} }{\hfill $\square$}
\newenvironment{pro*}[1][Proof]{{\it{#1:}} }{}

\newcounter{absatzcounter}[section]
\newenvironment{az}[1][]
{\refstepcounter{absatzcounter} \noindent {\bf (\thesection.\arabic{absatzcounter}) #1}}{\rm}

\begin{document}

\title{$\phantom{.}$ \\[-6ex] \hrule $\phantom{a}$ \\[3ex] \bf
      \Ueberschrift} 
  \author{Jakob Stix\thanks{\sc Mathematisches  Institut, 
   Universit\"at Bonn, Beringstra\ss e 1, 53115 Bonn  \newline
    \hspace*{0.45cm} E-mail address: {\tt stix@math.uni-bonn.de}
  \newline
  \noindent The author acknowledges support provided by DFG grant 
STI576/1-(1+2).
}}  
\date{\today}

\maketitle

\begin{quotation} 
  \noindent \small {\bf Abstract} --- Period and index of a curve $X/K$
 over a $p$-adic local field $K$ such that the fundamental group $\pi_1(X/K)$ admits a splitting are shown to be powers of $p$. As a consequence, examples of curves over number fields are constructed where having sections is obstructed locally at a $p$-adic place. Hence the section conjecture holds for these curves as there are neither sections nor rational points.
\end{quotation}


\section{Introduction} \label{sec:intro}

\begin{az}
This paper deals with the arithmetic of rational points of a smooth, projective curve $X$ over a field $k$ with algebraic closure $k^\alg$. By \textit{curve} we shall always mean a geometrically connected curve.
\end{az}

\begin{az} 
The \'etale fundamental group of a smooth, projective curve $X/k$ sits naturally in a short exact sequence 
\[ 1 \to \pi_1(X \otimes k^\alg) \to \pi_1(X) \to \Gal(k^\alg/k) \to 1,\]
which we abbreviate by $\pi_1(X/k)$. A $k$-rational point $x \in X(k)$ yields by functoriality a section of $\pi_1(X/k)$, with image the decomposition group of a point $\tilde{x}$ above $x$ in the universal pro-\'etale cover of $X$. Having neglected base points and due to the choice of $\tilde{x}$, only the class of a section up to conjugation by elements from $\pi_1(X \otimes k^\alg)$ is well defined. Let us denote by $S_{\pi_1{(X/k)}}$ the set of conjugacy classes of sections of $\pi_1(X/k)$. The section conjecture of Grothendieck's anabelian geometry \cite{letter}  speculates the following.
\begin{conj}[Grothendieck]
The natural map $X(k) \to S_{\pi_1(X/k)}$ which sends a rational point to the section given by its conjugacy class of decomposition groups is bijective if $k$ is a number field and the genus of $X$ is at least $2$.
\end{conj} 
There is also a version of the section conjecture for affine curves. Here rational points at infinity will lead to an abundance of additional \textit{cuspidal} sections. But apart from the obvious modification forced upon us by cuspidal sections the conjecture remains the same. The condition on the genus gets replaced by asking the Euler-characteristic to be negative. 

A birational version of the section conjecture over $p$-adic local fields was successfully addressed by Koenigsmann in \cite{Koenigsmann}.

Injectivity in the section conjecture is known \cite{Mz1} Theorem 19.1, \cite{EH} Theorem 5.1, see also \cite{Stix:pro-finite}, and only an exercise exploiting the Mordell-Weil Theorem. For surjectivity, the key claim, that is equivalent to the section conjecture, consists in \textit{having a section to be the only obstruction for  having a rational point}. Of course, this weaker assertion has to be available for all finite \'etale covers of a given curve in order to conclude.

The goal of the present paper is to provide evidence for the section conjecture. The evidence comes in form of theorems for which the assumption of $X$ having a rational point might be weakened to the extension  $\pi_1(X/k)$ being split. We will be able to provide examples of curves over number fields that have neither points nor sections of the fundamental group extension, see Section \ref{sec:trivial}. There will be a local obstruction that inhibits both sections and thus also rational points. The case of obstructions at real infinite places was already known \cite {Mz2} Theorem 3.13, see Section \ref{sec:trivial} and Appendix \ref{app:realSC}. 
\end{az}

\begin{az} 
The degree $\deg(x) = \dim_k(\kappa(x))$ of a closed point $x$ of a smooth projective curve $X/k$ is the degree of its residue field $\kappa(x)$ as an extension of $k$. The \textit{index} of $X$ is the greatest common divisor of the degrees of $k$-rational divisors:
\[ \ix(X) = \gcd(\deg(x); x \in X \text{closed}) = \#\Big(\bZ/\deg(\CH_0(X))\Big), \]
where the degree induces the degree map $\CH_0(X) \to \bZ$. The index equals $1$ if $X$ contains a $k$-rational point.
\end{az}

\begin{az}
The Picard variety classifying line bundles on $X$ decomposes as $\uPic_X = \coprod_{d \in\bZ} \uPic^d_X$ according to the degree $d$ of a line bundle. The part $\uPic_X^d$ of degree $d$ is a torsor under the abelian variety $\uPic_X^0$, that splits upon restriction to a field extension $k'/k$ if and only if $\uPic_X^d(k') \not= \emptyset$. Therefore the torsor is described by a class $[\uPic_X^d] = d \cdot [\uPic^1_X]$ in $\rH^1(k,\uPic^0_X)$. 

The \textit{period} of a smooth, projective curve $X/k$ is the order of the class $[\uPic^1_X]$, which is the smallest $d>0$ such that $\uPic_X^d(k) \not= \emptyset$ and which agrees with the greatest common divisor of the degrees of $k$-rational divisor classes of $X$. As for each closed point $x \in X$ we have 
\[\deg(x) \cdot [\uPic^1_X] = \cores_{\kappa(x)/k} [\uPic^1_{X \otimes \kappa(x)}] = 0,\]
the period divides the index. The period equals $1$ if $X$ contains a $k$-rational point.
\end{az}

\begin{az}
The \textit{period--index problem} raised by Lang and Tate \cite{LT} and studied further by Lichtenbaum \cite{Liindex} asks for conditions on $X$ or the ground field $k$ independent of the curve $X$ such that period and index agree. 
We will prove the following, see Theorem \ref{thm:powerofp}+\ref{thm:refined}. 
\begin{thm}
Let $K$ be a finite extension of $\bQ_p$ and let $X/K$ be a smooth, projective curve of positive genus, such that the fundamental group extension $\pi_1(X/K)$ admits a section. 

(1) For $p$ odd, $\pe(X)$ equals $\ix(X)$  and both are powers of $p$.

(2) For $p=2$, we have $\pe(X)$ and $\ix(X)$ are powers of $p$. If we 
moreover assume that we have an even degree finite \'etale cover $X \to X_o$ with $X_0$ of positive genus, then we have also $\pe(X)=\ix(X)$.   
\end{thm}
The following global application is Corollary \ref{cor:Q}.
\begin{thm}
Let $X$ be a smooth, projective curve over the rational numbers $\bQ$ of positive genus. If the fundamental group extension $\pi_1(X/\bQ)$ admits a section, then $\pe(X)$ equals $\ix(X)$.
\end{thm}
Of course, in light of the section conjecture, we expect $\pe(X)=\ix(X)$ equal to $1$ in both theorems above.
\end{az}


\section{The Brauer obstruction for line bundles} \label{sec:brauer}

Let $X/k$ be a smooth, projective variety over an arbitrary field $k$. We recall the theory of the Brauer obstruction for line bundles.

\begin{az}
A rational point $L \in \uPic_X(k)$ need not belong to a line bundle on $X$. After a finite Galois extension $k'/k$ the isomorphy class $L$ is realized by a line bundle $\cL \in Pic(X \otimes k')$ which moreover is $\Gal(k'/k)$ invariant. We find isomorphisms $\ph_\sigma : ^\sigma\cL \to \cL$ for every $\sigma \in \Gal(k'/k)$  where $^\sigma\cL = (\id \otimes \sigma)^\ast\cL$. The cocycle condition may be violated by the scalar automorphism
\[ (d\ph)_{\sigma,\tau} = \ph_\sigma  (^\sigma\ph_t) \ph_{\sigma\tau}^{-1} : \cL \to \cL \in  \rH^0\big(X\otimes k',\cAut(\cL)\big) = (k')^\ast \] 
not being the identity. Here the notation $d\ph$ is suggestive of the \v{C}ech boundary operator.  Of course $d(d\ph) =0$ and modifying the $\ph_\sigma$ corresponds to changing $d\ph$ by a coboundary. Moreover, the resulting class $b(L) = [d\ph] \in \rH^2(k,\Gm)$ does not depend on the choice of the splitting extension $k'/k$. The class $b(L)$ is the \textit{Brauer obstruction} that vanishes if and only if $\cL$ can be endowed with a descent datum relative $X \otimes k' \to X$ and hence if and only if $L$ belongs to a line bundle on $X$. The Brauer obstruction is additive for tensor products of line bundles, hence defines a homomorphism
$b: \uPic_X(k) \to \Br(k)$.
\end{az}

\begin{az}
The Leray spectral sequence in \'etale cohomology for the map $X \to \Spec(k)$ gives as exact sequence of low degree terms the following.
\[ 0 \to \Pic(X) \to \uPic_X(k) \xrightarrow{d_2^{0,1}} \Br(k) \to \Br(X) \]
The differential $d_2^{0,1}$ can be identified with the Brauer obstruction $b$ above, see Appendix \ref{app:ha}.
\end{az}

\begin{defi}
We define the \textbf{relative Brauer group} of a smooth, projective variety $X/k$ to be 
\[ \Br(X/k) = \ker\big(\Br(k) \to \Br(X)\big) \cong \uPic_X(k)/\Pic(X), \]
where the isomorphism is given by the Brauer obstruction map.
For smooth, projective curves $X/k$, the degree map defines a short exact sequence
\[ 0 \to \uPic_X^0(k)/\Pic^0(X) \to \Br(X/k) \xrightarrow{\deg} \pe(X)\bZ/\ix(X)\bZ \to 0.\]
We call the subgroup $\Br^0(X/k) = \uPic_X^0(k)/\Pic^0(X)$ of the relative Brauer group the \textbf{degree $0$ part} and the quotient $\ov{\Br}(X/k) = \pe(X)\bZ/\ix(X)\bZ$ the \textbf{N\'eron--Severi part} of the relative Brauer group.
\end{defi}
The period--index problem has an affirmative answer if and only if the N\'eron--Severi part of the relative Brauer group vanishes, e.g., if $\Br(k)$ vanishes.

\begin{az}
Without additional hypothesis on the base field $k$ we know that $\ix(X)$ divides $2g-2 = \deg(\Omega^1_{X/k})$, where $g$ is the genus of the curve $X$. Moreover, the index annihilates the relative Brauer group as for each closed point $x\in X$ the composite
\[ \Br(k) \to \Br(X) \xrightarrow{\ev_x} \Br(\kappa(x)) \xrightarrow{\cores_{\kappa(x)/k}} \Br(k)\]
is given by multiplication by $\deg(x)$. Here $\ev_x$ is the evaluation in the point $x \in X$.
\end{az}

\begin{az}
If the base field is finite, period and index are always $1$, although a curve over a finite field need not possess a rational point.
\begin{lem} \label{lem:pifinitefield}
Let $X$ be a smooth projective curve over a finite field $k$. Then $\ix(X)$ and $\pe(X)$ both equal $1$.
\end{lem}
\begin{pro}
Every torsor under an abelian variety over $k$ splits by a theorem of Lang, hence $\pe(X)$ equals $1$. As $\Br(k)$ vanishes, we must have $\ix(X) = \pe(X)$.
\end{pro}
\end{az}


\section{The geometry of Brauer obstructions} \label{sec:bs}

In this section we relate the Brauer obstruction for line bundles to the arithmetic  of Brauer--Severi varieties. For an introduction to Brauer--Severi varieties one may consult \cite{GS} Chapter 5, or \cite{ArtinBS}.

\begin{az}
A non-vanishing Brauer obstruction leads to `hyper'-structure on the line bundle.
\begin{lem} \label{lem:hyper}
Let $X/k$ be a smooth, projective variety, and let $L \in \uPic_X(k)$ be a Galois invariant line bundle. Let $\alpha \in \Br(k)$ be a Brauer class that is split by the Galois extension $k'/k$ and represented by the crossed product algebra $A$ corresponding to $\alpha \in \rH^2(k'/k,\Gm)$. The following are equivalent.
\begin{itemize}
\item[(a)] $b(L)= \alpha$.
\item[(b)] The line bundle $L$ is a genuine line bundle $\cL \in \Pic(X \otimes k')$, and $A$ acts semi-linearly on $q_\ast \cL$ where $q$ is the projection $q:X \otimes k' \to X$.
\end{itemize}
\end{lem}
\begin{pro} Let us assume (a). If $k'/k$ splits $\alpha$, then $b(q^\ast L) = 0$ and $L$ comes from a genuine line bundle $\cL$ on $X \otimes k'$. 
In the notation above, the not quite descent `cocycle' $\ph_\sigma$ for $L$ yields a $2$-cocycle $d\ph$ cohomologous to $\alpha$, so that we can write 
$A= \oplus_{\sigma \in \Gal(k'/k)} k' \cdot x_\sigma$  with 
\begin{eqnarray*}
x_\sigma x_\tau & = & (d\ph)_{\sigma,\tau} x_{\sigma\tau}\\
x_\sigma a & = & \sigma(a) x_\sigma \ \text{ for $a \in k'$.} 
\end{eqnarray*}
The algebra $A$ acts on $q_\ast \cL$ as follows. 
The scalars $k'$ act as scalars and $x_\sigma$ acts semilinearly via $q_\ast(\ph_\sigma) : q_\ast \cL = q_\ast (^\sigma \cL) \to q_\ast\cL.$ It is an easy computation, that these formulas define a $k'/k$-semilinear action of $A$ on $q_\ast \cL$.

On the other hand, if we assume $(b)$, then the action of $x_\sigma$ on $q_\ast \cL$ defines an isomorphism $\ph_\sigma : ^\sigma\cL \to \cL$, with which we compute the Brauer obstruction as $b(L) = [d\ph] = [A] = \alpha$.
\end{pro}
\end{az}

\begin{az}
Let $A$ be an Azumaya algebra of dimension $n^2$ over a field $k$. The associated Brauer--Severi variety $\BS_A$ over $k$ represents the functor on $k$-schemes
\[ \BS_A(T) = \{ A \otimes_k \OO_T \surj \cE \} \]
of simple $A \otimes_k \OO_T$-module quotients $\cE$, which are vector bundles of rank $n$, up to isomorphism.
The definition of a Brauer--Severi variety behaves well under base change, namely 
\[ \BS_{A \otimes_k k'} = \BS_A \otimes k'.\]
If the extension $k'/k$ splits $A$, then $\BS_A \otimes k'$ is the Brauer--Severi variety of the matrix ring $\Mat_n(k')$, which by Morita equivalence equals $\bP^{n-1}_{k'}$.  We choose to work with quotients as opposed to the more frequent left ideals in order to match Grothendieck's definition of projective space.
\end{az}

\begin{az}
Let $k'/k$ be a finite Galois extension of degree $m$ that splits $A$, and let $A \otimes \OO_{\BS_A} \surj \cE$ be the universal simple quotient. Let $q: \bP_{k'}^{n-1} = \BS_A \otimes k' \to \BS_A$ be the constant field extension that splits $\BS_A$. The sheaf $\OO(1)$ on $\bP_{k'}^{n-1} = \BS_A \otimes k'$ is Galois invariant and generates $\uPic_{\BS_A}(k)=\bZ$.
\begin{lem} \label{lem:bs}
We have $b(\OO(1)) = [A]$, the class of $A$ in $\Br(k)$. In particular $\Br(\BS_A/k)$ is the group $\langle [A] \rangle$ generated by $[A]$ in $\Br(k)$.
\end{lem}
\begin{pro}
This is stated in \cite{ArtinBS} Section 2 and attributed to Lichtenbaum in \cite{GS} Theorem 5.4.10. The proof only requires to compute $b(\OO(1))$, which we can easily accomplish using Lemma \ref{lem:hyper}, and which is therefore included here for the convenience of the reader.

By means of the standard idempotents $e_{i,i}$ in $\Mat_n(k') \cong A \otimes_k k'$ we may decompose the sheaf $q^\ast \cE$ into an $n$-fold direct sum of a line bundle $\cL$, see \cite{ArtinBS} 1.3. Morita equivalence translates the universal quotient on $\BS_{A \otimes k'}$ into $\OO^n \surj \cL$ for $\bP_{k'}^{n-1}$. Hence $\cL$ is nothing but $\OO(1)$. 

Let $k'\{G\}$ be the crossed product algebra for the trivial cocycle, hence $k'\{G\} \cong \Mat_m(k)$. On $q_\ast q^\ast \cE$ we have commuting $k$-linear actions of $A$ and $k'\{G\}$, hence an action of the algebra $A'= A \otimes_k k'\{G\}$. The idempotents $e_{i,i}$ lie in the subalgebra 
\[B = \Mat_n(k) \subseteq \Mat_n(k') \cong A \otimes_k k' \subseteq A'.\]
The centralizer $C=\rC_{A'}(B)$ of $B$ in $A'$ acts $k'/k$-semilinearly on $q_\ast(\OO(1)) = e_{i,i}(q_\ast q^\ast \cE)$. By the theorem of the centralizer, $C$ is an Azumaya algebra over $k$ of dimension $m^2$ and $[C] = [C] + [B] = [A'] = [A]$ in $\Br(k)$. The action of $C$ is exactly what Lemma \ref{lem:hyper} requires for the calculation of the Brauer obstruction: $b(\OO(1)) = [C] = [A]$.
\end{pro}
\end{az}

\begin{az}
The following result reveals the geometric origin of line bundles with Brauer obstruction.

\begin{prop} \label{prop:map}
Let $X/k$ be a smooth, projective curve, and let $A$ be an Azumaya algebra over $k$. The following are equivalent.
\begin{itemize}
\item[(a)] There is a map $f:X \to \BS_A$.
\item[(b)] There is a Galois invariant line bundle $L \in \uPic_X(k)$ with $b(L) = [A]$, in other words $[A] \in \Br(X/k)$.
\end{itemize}
\end{prop}
If $A$ is not a matrix algebra, the map $f$ in (a) must automatically be finite. Otherwise the image would be a $k$-rational point of $\BS_A$, saying that $A$ is split over $k$.

\begin{pro}
Let us assume (a). Then $L=f^\ast\OO(1)$ has $b(L) = b(\OO(1)) = [A]$ according to Lemma \ref{lem:bs}, so (b) follows.

Let us now assume (b). Let $k'$ be a maximal separable subfield of $A$. As $k'$ splits $A$, we have a genuine line bundle $\cL'$ on $X'=X \otimes k'$. Let $k"/k$ be a Galois extension that contains $k'$, and let $h:X" = X \otimes k" \to X'$, $q':X' \to X$ and $q" = q' \circ h$ be the projections. We apply Lemma \ref{lem:hyper} to the line bundle $\cL" = h^\ast \cL'$ and obtain an action of a $k"/k$-crossed product algebra $A(L) \cong \Mat_m(A)$ on $q"_\ast \cL"$ with $m=[k":k']$. 

The Galois group $H=\Gal(k"/k')$ acts on $q"_\ast\cL" = q'_\ast h_\ast h^\ast \cL'$ with invariants $q'_\ast \cL'$. Let $k"\{H\}$ be the trivial crossed product algebra over $k'$. It acts naturally on $k"$ so that $k"\{H\} = \End_{k'}(k") \cong \Mat_m(k')$. If we choose a normal basis for $k"/k'$ then the group $H$ sits inside the subalgebra $\Mat_m(k) \subseteq \Mat_m(k') \cong k"\{H\}$. By construction of $A(L)$, see Lemma \ref{lem:hyper}, $A(L)$ contains $k"\{H\}$ and thus $\Mat_m(k)$, such that $H \subset \Mat_m(k)$ acts as above. Hence the centralizer $C$ of $\Mat_m(k)$ in $A(L)$ acts on the $H$-invariants $q'_\ast \cL'$. By the theorem of the centralizer, $C$ is Brauer equivalent to $A$, and being of the same rank $C \cong A$. We have found an $A$ action on $\cE = q'_\ast \cL'$.

By computing ranks we deduce that the $A$ module $\cE$ on $X$ is simple and it remains to realize it as a quotient of $A \otimes \OO_X$. For that we may modify $L$ by actual line bundles from $\Pic(X)$. We may first assume, that $\rH^0(X,\cE) = \rH^0(X',\cL') \not= (0)$ contains a global section, which gives a non-trivial map $g:A \otimes \OO_X \to \cE$.
As fibrewise $\cE$ is a simple $A$-module, the map $g$ is either surjective or trivial in each fibre. The reduced locus $D=\{x ; g(x) \equiv 0\}$ is defined over $k$ and thus may be subtracted from $L$ leading to a new map 
\[g: A \otimes \OO_X \to \cE(-D)=q'_\ast \cL'(-D) \subseteq \cE.\] 
This process has to stop, because the degree of $L$ cannot become negative due to the presence of a global section. When it stops, we have found a simple quotient $A \otimes \OO_X \surj \cE$ that amounts to a map $f: X \to \BS_A$.
\end{pro} 

We see now again, that on a curve with a rational point no line bundle has a non-trivial Brauer obstruction. Indeed, by Proposition \ref{prop:map} we have a map $X \to \BS_A$ if $A$ is the value of a Brauer obstruction. The rational point on $X$ maps to a rational point on $\BS_A$, which shows that $[A]=0$ in $\Br(k)$.
\end{az}


\section{Period and Index over local fields} \label{sec:local}

In this section we begin to restrict the base field. We assume $X/K$ is a smooth, projective curve over a local $p$-adic field $K$, i.e, a finite extension $K/\bQ_p$. 

\begin{az}
The following result of Roquette has been given a proof by Lichtenbaum as simultaneously an intermediate step and a central ingredient of his duality for $p$-adic curves, see \cite{Liduality} Theorem 3.
\begin{thm}(Roquette; Lichtenbaum) \label{thm:roq}
Let $X/K$ be a smooth, projective curve over a finite extension $K/\bQ_p$. Then we have in $\Br(K) = \bQ/\bZ$ the following:
\[ \Br(X/K) = \bruch{1}{\ix(X)}\bZ/\bZ \text{\quad and \quad} \Br^0(X/K) = \bruch{1}{\pe(X)}\bZ/\bZ . \]
In other words, the relative Brauer group for curves over $p$-adic local fields is as big as it possibly could be.
\end{thm}
\end{az}

\begin{az}
Using Theorem \ref{thm:roq} Lichtenbaum, in \cite{Liduality} Theorem 7, provided the following set of conditions for the invariants period, index and genus of a curve over a $p$-adic local field. 
\begin{thm}(Lichtenbaum) \label{thm:Lichtenbaum}
Let $X/K$ be a smooth, projective curve of genus $g$ over a finite extension $K/\bQ_p$. Then we have  $\ix(X) | \pe(X) + 1 - g$, 
which implies

(1) $\pe(X) | g-1$,

(2) $\pe(X) | \ix(X) | 2\pe(X)$,

(3) if $\ix(X) = 2\pe(X)$ then $(g-1)/\pe(X)$ is odd.
\end{thm}
Theorem \ref{thm:Lichtenbaum} implies in particular an earlier result of Lichtenbaum's that the index of elliptic curves over $p$-adic local fields always equals the period, see \cite{Liindex} Theorem 3. 
\end{az}

\begin{az}
It was recently proven by Sharif \cite{Sharif} Theorem 2, that every triple $(g,\ix,\pe)$, which satisfies Lichtenbaum's conditions, occurs as the invariants genus, index, and period of a curve over a $p$-adic local field. Sharif uses branched, quadratic covers of torsors under elliptic curves for his examples.
\end{az}

\begin{az}
We would like to sketch the proof of Theorem \ref{thm:Lichtenbaum},
which at the beginning has nothing to do with $K$ being local $p$-adic. We therefore switch to a curve $X/k$. By Lemma \ref{lem:hyper}, for any 
$L \in \uPic_X(k)$ that is realized by a line bundle $\cL \in \Pic(X \otimes k')$, we have a semilinear Azumaya action on the sheaf $q_\ast\cL$, where $q:X \otimes k' \to X$ is the projection, by the $k'/k$-crossed product Azumaya algebra $A$ that represents $b(L)$. As this action subsists on cohomology, $\OO_X$-cohomology is as good as $A \otimes_k \OO_X$-cohomology, we find that
the order of $b(L)$ divides the index of $A$ which divides $\dim_{k'} \rH^i(X \otimes k',\cL)$ and moreover $\chi(X,L) = \deg(L) + 1-g$.

As the $L \in \uPic_X(k)$ of degree $\deg(L) = \pe(X)$ generate $\uPic_X(k)$ their images generate $\Br(X/k)$. The exponent of $\Br(X/k)$ therefore divides the lowest common multiples of the orders of such $b(L)$ and hence divides $\pe(X) + 1-g$. The rest follows from Theorem \ref{thm:roq} and uses crucially the fact that the base field is local $p$-adic.
\end{az}

\begin{az}
We would like to compute the index of a smooth, projective curve over a local $p$-adic field $K$ in terms of the special fibre of a good regular model. 

A result can be found in \cite{Clark} Theorem 9 in a setup, which is at the same time more general and more restrictive than ours, see also \cite{Clark} Section 4. The proof given below, which was obtained independently,  goes through unchanged to  include more general henselian discrete valuation rings with a perfect residue field.

\begin{thm} \label{thm:reduction}
Let $X \to S$ be a proper, flat curve with geometrically connected fibres over the spectrum $S$ of the ring of integers $\fo_K$ in a finite extension $K/\bQ_p$. We assume that the generic fibre $X_\eta$ is smooth over $K$, and that $X$ is regular. 

Let $Y=\bigcup_{\alpha \in A} Y_\alpha$ be the decomposition of the reduced special fibre $Y=X_{s,\redu}$ into irreducible components as a variety over $k$, the residue field of $\fo_K$ and set $Y_\alpha^0 = Y_\alpha - \Sing(Y)$ for the open part of $Y_\alpha$, where $Y$ is smooth over $k$, and set $q_\alpha : Y_\alpha^\norm \to Y_\alpha$ for the normalization map.  
Then we have 
\[ \ix(X_\eta) = 
\gcd_{\alpha \in A} (e_\alpha f_\alpha),\]
where $e_\alpha$ is the multiplicity of $Y_\alpha$ in the special fibre $X_s$ and $f_\alpha$ is the degree of $k_\alpha=\rH^0(Y_\alpha^\norm,\OO_{Y_\alpha^\norm})$ as an extension of $k$, which equals the field of constants of $Y_\alpha^0$.
\end{thm}
\begin{pro}
The closure of a point $x \in X$ is a horizontal divisor $T \subset X$, that meets the special fibre in only one point $y \in Y$. The map $T\to S$ being flat, we have 
\[\deg(x) = \deg(T/S) = (T \bullet X_s) = \sum_\alpha e_\alpha (T \bullet Y_\alpha) = \sum_\alpha e_\alpha  f_\alpha \deg_\alpha \big(q_\alpha^\ast \OO_X(T)\big),\]
where $(\ \bullet \ )$ is the intersection pairing in the sense of Lichtenbaum \cite{Licurves}, and $\deg_\alpha$ is the degree on $Y_\alpha^\norm$ with respect to the field of constants $k_\alpha$. This proves 
\[\gcd(e_\alpha f_\alpha; \alpha \in A) | \ix(X_\eta).\]

On the other hand, any point $y \in Y_\alpha^0$ can be lifted to be the intersection $\ov{\{x\}} \cap Y$ for a suitable $x \in X$, so that $T=\ov{\{x\}}$ is the zero locus of a regular parameter of $X$ in $y$, and thus $(T \bullet X_s) = e_\alpha f_\alpha \deg_\alpha(y)$.  We obtain that 
\[\ix(X_\eta) | \gcd(e_\alpha f_\alpha; \alpha \in A), \]
because the index of a curve over a finite field is $1$, see Lemma \ref{lem:pifinitefield}, and coprime multiples of the $0$-cycle of degree $1$ can be moved to avoid $Y^\norm_\alpha-Y_\alpha^0$ by Riemann-Roch. This proves the theorem.
\end{pro}
\end{az}


\section{The impact of a section on period and index} \label{sec:impact}

\begin{az}
The first observation that exploits the presence of a section works over general base fields.
\begin{prop} \label{prop:torsion}
Let $X/k$ be a smooth, projective curve such that $\pi_1(X/k)$ splits. Then the Brauer obstruction $b(L)$ vanishes for $L \in \uPic_X(k)$, which are torsion of order invertible in $k$.
\end{prop} 
\begin{pro}
We may assume that $L$ is of $\ell$-power order for some prime $\ell \not= \Char(k)$. The inclusion map $\bQ_\ell/\bZ_\ell(1) \inj \Gm$ leads to the following diagram with exact rows from the low degree terms of the Leray spectral sequence of $f:X \to \Spec(k)$.
\[\xymatrix@M=+1ex{ \rH^0\big(k,\rH^1(X\otimes k^\alg,\bQ_\ell/\bZ_\ell(1))\big) \ar[d]^{i} \ar[r]^(0.6){d_2^{0,1}} & \rH^2\big(k,\bQ_\ell/\bZ_\ell(1)\big) \ar[d] \ar[r]^{f^\ast} & \rH^2\big(X,\bQ_\ell/\bZ_\ell(1)\big) \ar[d] \\
\uPic_X(k) \ar[r]^b & \Br(k) \ar[r] & \Br(X) }
\]
It follows from the Kummer sequence $1 \to \mu_{\ell^n} \to \Gm \xrightarrow{\ell^n} \Gm \to 1$ in the limit over all $n$, that the image of $i$ equals the $\ell$-primary torsion $\uPic_X(k)[\ell^\infty]$. It remains to prove that $d_2^{0,1}=0$ or that $f^\ast$ is injective. But $f^\ast$ factors as
\[  \rH^2\big(k,\bQ_\ell/\bZ_\ell(1)\big) \xrightarrow{(\pi_1f)^\ast} \rH^2\big(\pi_1X,\bQ_\ell/\bZ_\ell(1)\big) \xrightarrow{c} \rH^2\big(X,\bQ_\ell/\bZ_\ell(1)\big),\]
with $(\pi_1f)^\ast$ being injective due to the section of $\pi_1(X/k)$. The comparison map $c$ is injective in general for $\rH^2$, which proves the proposition.
\end{pro}

\begin{rmk}
The above proof exploits the presence of a section in a very direct way. An alternative proof  is available by combining \cite{HS} Section 3.1 and \cite{Skoro} Corollary 2.39. 
\end{rmk}
\end{az}

\begin{az}
We return to the case of curves over local $p$-adic fields. 
\begin{cor} \label{cor:periodp}
Let $X/K$ be a smooth, projective curve over a finite extension $K$ of $\bQ_p$ such that $\pi_1(X/K)$ splits. Then $\pe(X)$ 
is a power of $p$.
\end{cor}
\begin{pro}
By Theorem \ref{thm:roq} we have to prove that the quotient $b: \uPic_X^o(K) \surj \Br^0(X/K)$ is a $p$-group. The group $\uPic_X^0(K)$ is isomorphic to $\fo_K^g \times \{\text{torsion}\}$ by a theorem of Mattuck \cite{Mattuck}. As $b$ kills torsion by Proposition \ref{prop:torsion} we obtain $\Br^0(X/K)$ as a quotient of a pro-$p$ group.
\end{pro}
\end{az}

\begin{az}
Having a section $s \in S_{\pi_1(X/K)}$ has another important impact. A section allows to construct geometrically connected finite \'etale covers. To an open subgroup $H \subset \pi_1(X)$ that contains $s(\Gal(k^\alg/k)$ is associated a geometrically connected cover $X_H \to X$ such that $\pi_1(X_H)=H$ and moreover the section $s$ lifts to a section $s_H$ of $H \surj \Gal(k^\alg/k)$.  The construction of $X_H$ involves choices of base points. 

We call such a finite \'etale cover that admits a lifting of the section $s$ a \textit{neighbourhood} of $s$. There is an abundance of neighbourhoods of sections as every closed subgroup of a pro-finite group is the intersection of all open subgroups that contain it. Using covers of even degree we can prove the following.
\begin{thm} \label{thm:powerofp}
Let $X/K$ be a smooth, projective curve of genus $g>0$ over a finite extension $K$ of $\bQ_p$ such that $\pi_1(X/K)$ splits. Then $\ix(X)  = \#\Br(X/K)$
is a power of $p$.
\end{thm}
\begin{pro}
We need to show that $\Br(X/K)$ is a $p$-group. As for a finite \'etale cover $X' \to X$ we have $\Br(X/K) \subseteq \Br(X'/K)$ we may replace $X$ by a neighbourhood of a section of $\pi_1(X/K)$. If we choose $\deg(X'/X)$ to be even, the genus $g'$ of $X'$ will become odd by the Hurwitz formula $g'-1 = \deg(X'/X) (g-1)$. Such covers exist by the assumption on the genus. We may and will therefore assume that the genus is odd.
By Corollary \ref{cor:periodp} and Theorem \ref{thm:Lichtenbaum}(2) it suffices to treat the case $p \not= 2$. Thus $(g-1)/\pe(X)$ is even and by Theorem \ref{thm:Lichtenbaum}(3) we must have $\ix(X)=\pe(X)$ which is a power of $p$ by Corollary \ref{cor:periodp}.
\end{pro}
\end{az}

\begin{az}
The above proof actually shows more.
\begin{thm} \label{thm:refined}
Let $X/K$ be a smooth, projective curve of genus $g>0$ over a finite extension $K$ of $\bQ_p$ such that $\pi_1(X/K)$ splits. 

(1) If $p$ is odd, we have $\pe(X)=\ix(X)$.

(2) If $p=2$, we have $\pe(X')=\ix(X')$ for every finite \'etale cover $X' \to X$ of even degree such that $\pi_1(X'/K)$ splits.
\end{thm}
\begin{pro}
The case of odd $p$ follows from Theorem \ref{thm:Lichtenbaum}(2) as both period and index are powers of $p$. We therefore assume $p=2$ in the sequel.

\textit{Step 1:} We have $\pe(X)=\ix(X)$ if and only if $\ix(X) | g-1$. 

\noindent Indeed, $\pe(X)$ always divides $g-1$ by Theorem \ref{thm:Lichtenbaum}(1)and so does $\ix(X)$ if period and index agree. On the other hand, if $\ix(X)|g-1$ then $\ix(X)/\pe(X)$ divides $(g-1)/\pe(X)$ and Theorem \ref{thm:Lichtenbaum}(3) only allows $\pe(X)=\ix(X)$.

\textit{Step 2:} Next we discuss obstructions for square roots. From the multiplication by $2$ sequence of $\uPic^0_X$ we get a short exact sequence in Galois cohomology
\[ 0 \to \uPic_X^0(K)/2\uPic_X^0(K) \to \rH^1(K,\uPic_X^0[2]) \to \rH^1(K,\uPic^0_X)[2] \to 0.\]
For $M \in \Pic^{2g-2}(X)$ as an element of $\uPic_X^{2g-2}(K)$ we get via the boundary map in Galois cohomology of the sequence $0 \to \uPic_X^0[2] \to \uPic_X \xrightarrow{2\cdot} \uPic_X^{\rm even} \to 0$ an obstruction $[\sqrt{M}] \in \rH^1(K,\uPic_X^0[2])$ that vanishes if and only if $M$ has a Galois invariant square root. From the commutativity of the following diagram
\[\xymatrix{ 0 \ar[r] &  \uPic_X^0[2] \ar[d] \ar[r] & \uPic_X \ar[d] \ar[r]^{2 \cdot} &  \uPic_X^{\rm even} \ar[d]^{\bruch{1}{2} \deg} \ar[r] & 0 \\
 0 \ar[r] &  \uPic_X^0  \ar[r] & \uPic_X \ar[r] &  \bZ  \ar[r] & 0 
}\]
we conclude that $[\sqrt{M}]$ maps to $[\uPic_X^{g-1}]$ in $\rH^1(K,\uPic_X^0)[2]$ which vanishes because $\pe(X)$ divides $g-1$. Therefore the set of all $[\sqrt{M}]$ forms a coset in $\uPic_X^0(K)/2\uPic_X^0(K)$ under the image of $\Pic^0(X)$ under the canonical map. The coset is nonempty as it contains the square root obstruction to the canonical bundle, that is the obstruction to a Galois invariant spin structure or $\theta$-characteristic. We arrive at a well defined class $\rho = \rho_{X/K}$ in the cokernel $\Br^0(X/k) \otimes \bZ/2\bZ$ that vanishes if and only if a suitable choice of $M \in \Pic^{2g-2}(X)$ has a square root in $\uPic_X^{g-1}(K)$. The formation of the class $\rho$ clearly is compatible with pullback under finite \'etale covers $f:X' \to X$, where $f^\ast\rho_{X/K} = \rho_{X'/K}$.

\textit{Step 3:} If $\rho = 0$ and $g-1$ is even, then $\ix(X) | g-1$.

\noindent Indeed, if $\rho$ vanishes then we have $L \in \uPic_X^{g-1}(K)$ such that $2L$ is a line bundle on $X$. The Brauer obstruction of $L$ satisfies $2b(L) = b(2L) = 0$ and thus is $0$ or $\bruch{1}{2}$. If $b(L) = 0$, then $L$ is a line bundle and $\ix(X)$ divides $g-1$.

If on the other hand $b(L)= \bruch{1}{2}$ and $\pe(X) \not= 1$, then $\bruch{1}{2} \in \Br^0(X/K)$ and we find $L_0 \in \uPic^0_X(K)$ with $b(L) = b(L_0)$. But then $L-L_0$ has degree $g-1$ and is a line bundle on $X$ as $b(L-L_0) = 0$, so that $\ix(X)$ again divides $g-1$. 

It remains to deal with $\pe(X) = 1$. But then $\ix(X)$ divides $2$ which divides $g-1$ by assumption.

\textit{Step 4:} For a finite \'etale cover $f:X' \to X$ of even degree, the Hurwitz formula $g'-1 =\deg(X'/X) (g-1)$ for the genus $g'$ of $X'$ gives the first condition of step 3. By step 1 and 3, it suffices to see $\rho_{X'/K}$ vanishing. If the pullback map $\Br^0(X/K) \inj \Br^0(X'/K)$ is not surjective, then $f^\ast \rho_{X/K}$ is divisible by $2$ and $\rho_{X'/K}$ vanishes. Otherwise we have $\pe(X) = \pe(X')$ and hence 
\[\ix(X') | 2\pe(X') = 2\pe(X) | \deg(X'/X) (g-1) = g'-1\]
and we conclude with step 1.
\end{pro}

We remark that when $p=2$, ramified covers $X \to \BS_{1/2}$ of odd degree, where $\BS_{1/2}$ is the Bauer--Severi variety of the skew field with invariant $\bruch{1}{2}$, lead to curves, where $\rho_{X/K} = 0$, $\pe(X)$ is odd and $\ix(X) = 2\pe(X)$. The section conjecture for $p$-adic base fields prohibits that $\pi_1(X/K)$ splits for such curves, but we are unable to prove this at the moment.
\end{az}


\section{Global Brauer obstructions} \label{sec:global}

\begin{az}
In characteristic $0$ or for proper varieties, the fundamental group exact sequence can be base changed due to the geometricity of the fundamental group. In other words, for an extension $E/F$ of fields and a smooth, projective curve $X/F$ the pullback $f^\ast\pi_1(X/F)$ along the natural restriction map $f: \Gal(E^\alg/E) \to \Gal(F^\alg/F)$ is canonically isomorphic to the extension $\pi_1(X \otimes E/E)$. In particular, we get a base change map of spaces of sections:
\[ S_{\pi_1(X/F)} \to S_{\pi_1(X \otimes E/E)}, \qquad s \mapsto s_E.\]
\end{az}

\begin{az}
In this section we work with curves $X$ over number fields, which are called $F$ for distinction purposes. The base change of a section $s$ of $\pi_1(X/F)$ to the completion $F_\fp$ of $F$ with respect to a place $\fp$ gives a section $s_\fp = s_{F_\fp}$ of $\pi_1(X \otimes F_\fp/F_\fp)$. Hence we may exploit the local results that were obtained in Section \ref{sec:impact}. 
\begin{thm} \label{thm:global}
Let $X/F$ be a smooth, projective curve over a number field $F$ such that 
\begin{enumerate}
\item $\pi_1(X/F)$ admits a section, and
\item for each prime number $p$ there is at most one places $\fp$ of $F$ with residue characteristic $p$ such that $X$ has bad reduction at $\fp$.
\end{enumerate}
Then there is no Brauer obstruction for line bundles: $\Br(X/F) = 0$, in other words $\Pic(X) = \uPic_X(F)$ and the map $\Br(F) \to \Br(X)$ is injective. In particular, $\pe(X)$ equals $\ix(X)$.
\end{thm}
\begin{pro} Let $\fp$ be a finite place of $F$. 
If $X_\fp = X \otimes F_\fp$ has good reduction, we have $\ix(X_\fp)=1$ and $\Br(X_\fp/F_\fp) = 0$ by Theorem \ref{thm:reduction}. Regardless of the reduction behaviour, by base change $\pi_1(X_\fp/F_\fp)$ admits a section and $\Br(X_\fp/F_\fp)$ is a $p$-group by Theorem \ref{thm:powerofp} where $p$ is the residue characteristic at $\fp$. 

If $\fp$ is a real infinite place, the section $s_\fp$ forces $X_\fp$ to have a rational point by the real section conjecture, see Theorem \ref{thm:sip}. So for infinite places $\fp$ the group $\Br(X_\fp/F_\fp)$  vanishes.

The restriction of a class $\alpha \in \Br(X/F)$ to a completion $F_\fp$ lies in $\Br(X_\fp/F_\fp)$. Thus the local invariants $\inv_\fp(\alpha)$ are of prime power order with each prime occuring at most once. The Hasse-Brauer-Noether local-global principle for Brauer groups requires the sum of the local invariants to vanish, which is impossible unless $\alpha=0$. This proves the theorem.
\end{pro}

\begin{cor} \label{cor:Q}
On a smooth, projective curve $X$ over $\bQ$, the existence of a section of $\pi_1(X/\bQ)$ implies that every Brauer obstruction for line bundles vanishes. In particular, $\pe(X)$ equals $\ix(X)$. 
\end{cor}
\begin{pro}
The second condition of Theorem \ref{thm:global} is automatically satisfied for $F=\bQ$.
\end{pro}
\end{az}




\section{Trivial examples for the section conjecture} \label{sec:trivial}

\begin{az}
An example for the section conjecture is a smooth, projective curve $X$ of genus $\ge 2$ over a number field $F$, such that every section of $\pi_1(X/F)$ belongs to the decomposition group of an $F$-rational point of $X$. Our examples will be trivial in the sense, that the curves we provide will have neither sections nor rational points. Although this seems to be `empty mathematics', understanding 
the case that has no point is the key to settle the section conjecture.
\end{az}

\begin{az}
The example we offer, will always be obstructed locally at either a finite or an infinite place. So the reason for not having a section is purely local. As is well known, there are curves over number fields, such that everywhere locally we have points but nevertheless there is no global point. This illustrates the limitations of our method to construct examples, which otherwise is fairly complete.

Another limitation consists in not varying the base field $F$ for a fixed curve $X/F$. The base changes $X \otimes F'$ will attain rational points eventually and leave the realm of our examples.
\end{az}

\begin{az}
Local obstructions at infinite real places are due to the real section conjecture, see Appendix \ref{app:realSC}. By Theorem \ref{thm:sip}, a curve $X/F$ such that $X(\bR) = \emptyset$ for some real infinite place $F \inj \bR$ cannot have sections of its $\pi_1(X/F)$ and thus the section conjecture holds trivially. 
\begin{ex}
Any branched cover $X$ of genus $\ge 2$ of the curve $\{X^2+Y^2+Z^2 = 0\} \subset \bP^2_F$ over a number field $F$ with a real place, e.g., 
\[\{X^{2n} + Y^{2n} + Z^{2n} = 0\} \]
for $n \geq 2$ and $F=\bQ$, satisfies the section conjecture trivially. 
\end{ex}
It is not a coincidence that a Brauer--Severi variety appears in the above example, see Corollary \ref{cor:norealscover}. The same will happen for the local $p$-adic examples that we construct.
\end{az}

\begin{az}
Let $F$ be a number field with a finite place $\fp$ of residue characteristic $p$. Let $A$ be an Azumaya algebra over $F$ with $\inv_\fp(A)$ of order not a power of $p$. Let $X/F$ be a smooth, projective curve with a map $X \to \BS_A$. Such a curve may be constructed by Bertini's Theorem by suitable general hyperplane sections of $\BS_A \subseteq \bP^N_F$ of high degree. By base change, we obtain a map $X_\fp = X \otimes F_\fp \to \BS_{A \otimes F_\fp}$, so that $\Br(X_\fp/F_\fp)$ contains $[A \otimes F_\fp]$ by Proposition \ref{prop:map}, which is an element of order not a power of $p$ by construction. Therefore Theorem \ref{thm:powerofp} prevents $\pi_1(X_\fp/F_\fp)$ and thus also  $\pi_1(X/F)$ from admitting a section. The section conjecture holds for $X/F$ trivially, being obstructed locally at the finite place $F \inj F_\fp$.

\begin{ex}
Let $p$ be odd and let $a,b \in \bQ^\ast$ be not both negative with non trivial Hilbert symbol $(a,b)_p$ locally at $\bQ_p$. For $p \equiv 3 \mod 4$ we can take for example $(a,b) = (p,-1)$. The curve  $\BS_{a,b} = \{X^2 - aY^2 -bZ^2 = 0\} \subseteq \bP^2_{\bQ}$ gives upon base change to $\bQ_p$ the Brauer--Severi variety for the Quaternion algebra with invariant $\bruch{1}{2}$. Any branched cover $f:X \to \BS_{a,b}$ with genus $\geq 2$, e.g.,
\[ \{X^{2n} - aY^{2n} - bZ^{2n} = 0\}, \]
as a curve over $\bQ$, verifies the section conjecture trivially, being obstructed locally at $p$. If the degree of $f$ is odd, the assumption on the sign of $a,b$ assures that we have real points, so that these examples are not accounted for by obstructions at a real infinite place.
\end{ex}
\end{az}

\begin{az}
Let us provide an alternative construction of trivial examples for the section conjecture that is based on the computation of the index from the degeneration of a curve.

Let $k$ be a finite field of characteristic $p$ and let $C/k$ be a curve without automorphisms:   $\Aut(C \otimes k^\alg)= 1$. Let $k_n/k$ be the unique degree $n$ extension and let $\Frob$ be the Frobenius in $\Gal(k_n/k)$. Take $y_1,y_2 \in C(k_n)$ not conjugate to each other such that the residue field of $C$ at $y_i$ is $k_n$. Then glue $C \otimes k_n$ at the corresponding $k_n$-rational points $y_1,y_2$ in a transversal manner using $\Frob$ to identify the residue fields at $y_1$ and $y_2$. For the result of the gluing $Y=C \otimes k_n/(y_1 \sim y_2)$ a neighbourhood of the identified point $y$ is given by an open  $U \subseteq  C \otimes k_n$ such that both $y_i \in U$ with sheaf of functions 
\[\OO_Y(U) = \{f \in \OO_C(U) ; f(y_1) = \Frob(f(y_2))\}.\]
The curve $Y$ is semistable over $k$ with normalization $C \otimes k_n$. Geometrically, we have $n$ conjugate copies of $C$ that are glued transversally along the conjugates of the points $y_1$ with the next conjugate of $y_2$, so that the dual graph is a circle of length $n$.

The curve $Y/k$ has only one component with $f=n$ and being reduced $e=1$. We pick $n$ to be not a power of $p$ and deform $Y$ to a smooth curve $X$ over a 
number field $F$, such that $k$ is the residue field at a place $\fp$ and the fibre of a regular model above $\fp$ is $Y$. This can be done as  $\ov{\mg}$ is smooth over $\bZ$ with $\mg \subset \ov{\mg}$ dense in each fibre. Taking care of $Y$ not having automorphisms allows to deform without enlarging the residue field. By Theorem \ref{thm:reduction} the index of $X_\fp$ is $n$, which we chose not to be a power of $p$. 
This is the strategy of \cite{Clark} for providing curves over local fields with interesting index. In our case, we conclude that $\pi_1(X/F)$ does not admit sections with again an obstruction provided locally at $\fp$. Hence the curve $X/F$ satisfies the section conjecture trivially.
\end{az}


\begin{appendix}

\section{The real section conjecture} \label{app:realSC}

\begin{az}
We discuss the section conjecture for smooth, projective curves $X/\bR$ over the reals $\bR$ for matters of completeness and also to show the technique of the main body of the paper in this case. 

The real section conjecture was formulated and proven already by Mochizuki in \cite {Mz2} Theorem 3.13. Huisman in \cite{Huisman}  Theorem 7.2 determined the structure of the fundamental group without a point, which as well as a result of Cox, \cite{Cox} Theorem 2.1, going back to Artin--Verdier from 1964, could serve as the key ingredient of surjectivity. We will rely on a theorem of Witt and the discussion of the relative Brauer group.
\end{az}

\begin{az}
The most difficult part in the section conjecture is to recognize the presence of rational points. For real curves we have the following criterion.
\begin{thm} \label{thm:sip}
Let $X/\bR$ be a smooth, projective curve of genus $>0$. Then the following are equivalent.
\begin{itemize}
\item[(a)] $\pi_1(X/\bR)$ admits a section.
\item[(b)] $\Br(X/\bR) = (0)$.
\item[(c)] $X(\bR) \not= \emptyset$.
\end{itemize}
\end{thm}
\begin{pro}
That (b) implies (c) follows from the following theorem of Witt, see \cite{Scheiderer} 20.1.3.
\begin{thm}[Witt 1934]
For a smooth, projective curve $X/\bR$ the map 
\[ \Br(X) \to \rH^0\big(X(\bR), \bZ/2\bZ\big) \]
given by pointwise evaluation is an isomorphism.
\end{thm}
Namely, if $\Br(X/\bR)$ vanishes, then $\Br(X)$ contains the constant classes $\Br(\bR) \inj \Br(X)$ and hence $X(\bR)$ cannot be empty.

That (c) implies (a) is obvious. So let us assume (a) and prove (b). As real connected Lie groups are divisible, we find 
\[\uPic_X^0(\bR) = \Big(\uPic_X^0(\bR)\Big)_0 \times \pi_0\big(\uPic_X^0(\bR)\big)\]
where $\Big(\uPic_X^0(\bR)\Big)_0$ is the connected component of $\uPic_X^0(\bR)$ and $\pi_0\big(\uPic_X^0(\bR)\big)$ is the finite group of connected components. 
By Proposition \ref{prop:torsion} the Brauer obstruction vanishes on torsion line bundles, so the $\Pic^0$-part of $\Br(X/\bR)$ vanishes and 
\[\Br(X/\bR) = \ov{\Br(X/\bR)} = \pe(X)\bZ/\ix(X)\bZ.\]
Let $f:X' \to X$ be a neighbourhood of a section of $\pi_1(X/\bR)$ of even degree $\deg(f)$. Then $f^\ast:\Br(X/\bR) \to \Br(X'/\bR)$ is injective and given by multiplication by the degree: 
\[\deg(f) \cdot : \pe(X)\bZ/\ix(X)\bZ \to \pe(X')\bZ/\ix(X')\bZ.\]
As $\ix(X')$ is at most $2$ the map $f^\ast$ is the zero map showing that $\Br(X/\bR)$ vanishes. 
\end{pro}

\begin{cor} \label{cor:norealscover}
A smooth, projective curve $X/\bR$ without real points is a finite branched cover of $\BS_\bH = \{X^2 + Y^2 + Z^2 = 0\} \subset \bP^2_\bR$.
\end{cor}
\begin{pro}
Indeed, if $X(\bR) = \emptyset$ then by Theorem \ref{thm:sip} we find $[\bH] \in \Br(X/\bR)$ from which Proposition \ref{prop:map} constructs a cover of $\BS_\bH$.
\end{pro}

As an alternative, one may refine \cite{Cox} Theorem 2.1 to the following Theorem.
\begin{thm}
Let $X/\bR$ be a smooth, projective curve of genus $>0$. Then the following are equivalent.
\begin{itemize}
\item[(a)] $\pi_1(X/\bR)$ admits a section.
\item[(b)] $\pi_1(X)$ contains non-trivial $2$-torsion.
\item[(c)] The $2$ cohomological dimension of $\pi_1(X)$ is unbounded.
\item[(d)] $X(\bR) \not= \emptyset$.
\end{itemize}
\end{thm}
\begin{pro}
The image of a section is $2$-torsion, and $2$-torsion forces the $2$-cohomological dimension to be unbounded, so (a) $\Longrightarrow$ (b) $\Longrightarrow$ (c). As (d) obviously implies (a), it remains to assume (c) and to prove (d). Since the cohomology of the fundamental group of a curve of positive genus coincides with the \'etale cohomology of the curve, this is an immediate consequence of the calculation of the \'etale cohomology of a real curve in terms of equivariant cohomology as done by Cox in the following theorem.
\end{pro}

\begin{thm}[Cox, \cite{Cox} Proposition 1.2, see \cite{Scheiderer} Introduction]
Let $X/\bR$ be a smooth, projective variety of dimension $d$. Then for $q>2d$ we have
\[ \rH^q(X,\bZ/2\bZ) \cong \bigoplus_{i=0}^d \rH^i\big(X(\bR),\bZ/2\bZ\big).\] 
\end{thm}
\end{az}

\begin{az} 
We recall that a neighbourhood of a section $s \in S_{\pi_1(X/k)}$ consists of a connected finite \'etale cover $X' \to X$ such that $s(\Gal(k^\alg/k))$ is contained in $\pi_1(X') \subset \pi_1(X)$. More precisely, we are forced to work with pointed covers, but we ignore this detail here. The system of all neighbourhoods of $s$ forms the decomposition tower $X_s$ of the section $s$ (not the class of a section), the projective system $X_s = (X')$ of all connected finite \'etale covers $X' \to X$ together with a lift of the section. It makes sense to speak of the fundamental group of this pro-\'etale cover. We find $\pi_1(X_s) = \bigcap_{X'} \pi_1(X') = s\big(\Gal(k^\alg/k)\big)$, where $X'$ runs through the system of all neighbourhoods of $s$. 

It now follows formally, that two conjugacy classes of sections agree if and only if they have the same
neigbourhoods. We could afford to be sloppy with the distinction between a section and a class of sections, as well as with covers and pointed covers, because the conjugacy action only moves the corresponding chosen lift of a section but does not affect the possibility to lift a section to a neighbourhood.

For a geometric section $s_x$ that belongs to a rational point $x \in X(k)$ a neighbourhood is a cover $X' \to X$ such that $x$ lifts to a $k$-rational point $x' \in X'(k)$. 

For surjectivity of the section conjecture we need to find a coherent system of lifts of a rational point in the decomposition tower $X_s$ of the section $s$.

For injectivity, we need to find for each pair $x,y$ or rational points a finite \'etale cover $X'\to X$ such that exactly one of the points lifts. 
\end{az}

\begin{az}
Let $X/\bR$ be a smooth, projective real curve. A finite \'etale map $X' \to X$ induces a proper topological covering map $X'(\bR) \to X(\bR)$. So if $x \in X(\bR)$ lifts to $X'(\bR)$ then the whole connected component of $X(\bR)$ that contains $x$ does. Thus the neighbourhoods of two points in the same connected components agree, and the map of the section conjecture factors over the set of connected components
\begin{thm}[Real section conjecture]
Let $X/\bR$ be a smooth, projective curve of genus $\geq 1$. Then the map 
\[ \pi_0\big(X(\bR)\big) \to S_{\pi_1(X/\bR)},\]
that maps a connected component of the real locus $X(\bR)$ to the corresponding conjugacy class of sections is a bijection of finite sets.
\end{thm}
\begin{pro}
We have already seen, that the map is well defined.
For a section $s \in S_{\pi_1(X/\bR)}$ to be in the image, we need to find a rational point in $\varprojlim X'(\bR)$, where the limit extends over the projective system of all neighbourhoods of the section $s$. By Theorem \ref{thm:sip} all these $X'(\bR)$ are nonempty. The limit is nonempty as a limit of nonempty compact sets. This proves surjectivity.

In order to prove injectivity, we consider the canonical embeding $X \inj \uPic_X^1$ and the Kummer sequence for $\uPic_X^0$ as follows.
\begin{lem} \label{lem:maptopic1}
The map $X(\bR) \to \uPic_X^1(\bR)/2\uPic_X^0(\bR)$ of sets factors over an injective map
\[\pi_0\big(X(\bR)\big) \inj \uPic_X^1(\bR)/2\uPic_X^0(\bR).\]
\end{lem}
\begin{pro}[Proof of the lemma]
Let $x,y \in X(\bR)$ map to the same $2\uPic^0_X(\bR)$-coset of $\uPic_X^1(\bR)$. Then $x-y \sim 2D$ with a divisor $D$ in a Galois invariant divisor class. As here $X$ has real points, all Brauer obstructions vanish and we may choose $D$ to be Galois invariant itself. So there is a function $f \in \bC(X)^\ast$ with Galois invariant $\divisor(f) = x-y-2D$, which implies that $c_\sigma=\sigma(f)/f \in \bC^\ast$ for every $\sigma \in \Gal(\bC/\bR)$. The corresponding $1$-cocycle $c_\sigma$ of $\Gal(\bC/\bR)$ with values in $\Gm$ is a coboundary by Hilbert's Theorem 90. We may modify $f$ by a constant accordingly, so that $\sigma(f) = f$, or $f \in \bR(X)^\ast$, which means that $f$ takes real values at real points of $X$. 

The function $f$ when restricted to $X(\bR)$ is continuous real valued with sign changes exactly at poles or zeros of odd order, so exactly in $x$ and $y$. It follows that $x$ and $y$ lie on the same connected component.

That the map of the lemma indeed factors over the set of connected components $\pi_0(X(\bR))$ follows from the fact, that the section of a point is constant on connected components and the discussion of the Kummer sequence below.
\end{pro}

We continue our proof of injectivity.
By Lemma \ref{lem:maptopic1} and a translation argument we may replace $X$ by the abelian variety $A=\uPic^0_X$. It remains to show that any $a \in A(\bR)$ with section $s_a$ conjugate to the section $s_0$ associated to the origin actually lies in $2A(\bR)$. 
The Kummer sequence $0 \to A[2] \to A \to A \to 0$
induces an injective boundary map 
\[\delta : A(\bR)/2A(\bR) \inj \rH^1(\bR,A[2]).\]
The proof of injectivity will therefore be completed by the following Proposition.
\end{pro}
\end{az}

\begin{az}
Let $A/k$ be an abelian variety and $n$ invertible in $k$. We consider the Kummer sequence $0 \to A[n] \to A \to A \to 0$ and the truncated mod $n$ fundamental group extension $\pi_1^{[n]}(A/k)$
\[ 1 \to A[n] \to \pi_1^{[n]}A \to \Gal(k^\alg/k) \to 1,\]
which is a quotient of $\pi_1(A/k)$. Hence rational points $a\in A(k)$ lead to sections of $\pi_1(A/k)$ which project onto sections $s_a$ of $\pi_1^{[n]}(A/k)$.
\begin{prop}
The boundary of the Kummersequence 
\[\delta : A(k)/nA(k) \inj \rH^1(k,A[n])\]
maps $a \in A(k)$ to the class of the difference cocycle $\sigma \mapsto s_a(\sigma) s_0(\sigma)^{-1}$ of the section $s_a$ and the section $s_0$ associated to the origin.
\end{prop}
\begin{pro}
It is sufficient to compute and compare the actions of $\Gal(k^\alg/k)$ on the fibre $A[n](k^\alg)$ of multiplication by $n$ over $0$ that are induced by $s_a$ and $s_0$ respectively. Via $s_0$ the action is the natural one.

The same holds for $s_a$ for the action on the fibre over $a$. In order to regard $s_a$ as a section of $\pi_1(A,0)$, so we can let it act on the fibre above $0$, we have to choose a path from $0$ to $a$. The path is determined by what it does for the cover `multiplication by $n$' and thus is given by the translation by a fixed $n$th root $\bruch{1}{n}a$ of $a$. 

For a point $P \in A[n](k^\alg)$ we compute
\[ s_a(\sigma)\big(P\big) = - \bruch{1}{n}a + \sigma\big(P + \bruch{1}{n}a\big) = - \bruch{1}{n}a + \sigma\big(P\big) + \sigma\big(\bruch{1}{n}a\big) =  \sigma\big(\bruch{1}{n}a\big) - \bruch{1}{n}a + s_0(\sigma)\big(P\big).\]
Hence $s_a(\sigma)s_0(\sigma)^{-1}$ equals translation by $\sigma(\bruch{1}{n}a) - \bruch{1}{n}a$, which is nothing but the value of a cocycle for $\delta(a)$ at $\sigma$.
\end{pro}
\end{az}


\section{Brauer obstruction and Leray spectral sequence} \label{app:ha}

\begin{az}
In this appendix we compute the differential $d_2^{0,1}$ in the Leray spectral sequence of a smooth, projective variety $f:X \to \Spec(k)$ for coefficients $\Gm$. This is the map 
\[ d_2^{0,1} : \uPic_X(k) = \rH^0\big(k,\rH^1(X \otimes k^\alg,\Gm)\big) \to \Br(k) = \rH^2\big(k,\rH^0(X \otimes k^\alg,\Gm)\big).
\]
\begin{prop}
We have $b(L) = d_2^{0,1}(L)$ for each $L \in \uPic_X(k)$, where $b(L)$ is the Brauer obstruction for line bundles.
\end{prop}
\end{az}

\begin{az}
The Leray spectral sequence is a Grothendieck spectral sequence for a composite $F \circ G$ of functors, where $F= \rH^0(k,-)$ and $G = f_\ast$. It comes from a natural filtration on the total complex of $F J^{\bullet, \bullet}$, where $J^{\bullet,\bullet}$ is an Eilenberg--MacLane resolution $GI^\bullet \to J^{\bullet,\bullet}$, the resolution $M \to I^\bullet$ computes the higher direct images $R^qG(M)$, and $M$ are the coefficients we are interested in. 

The double complex has two differentials $d'$ and $d"$ which anticommute, so that $d=d'+d"$ furnishes the total complex with a differential. Each group in the spectral sequence is a subquotient of the total complex $\Tot(FJ^{\bullet,\bullet})$ and the differential of the spectral sequence is induced by $d$. 
\[\xymatrix@R=1pc@C=1pc{ \vdots & \vdots & \vdots &  \\
FJ^{0,1} \ar[r]^{d'} \ar[u]^{d"} & FJ^{1,1} \ar[r]^{d'} \ar[u]^{d"} & FJ^{2,1}\ar[r]^{d'} \ar[u]^{d"} & \cdots \\
FJ^{0,0} \ar[r]^{d'} \ar[u]^{d"} & FJ^{1,0} \ar[r]^{d'} \ar[u]^{d"} & FJ^{2,0} 
\ar[r]^{d'} \ar[u]^{d"} & \cdots 
} \qquad
\xymatrix@R=1pc@C=1.5pc{  0  &  & &  \\
x \ar@{|->}[r]^(0.42){d'} \ar@{|->}[u]^{d"} & d'(x)  &  & \\
 & y \ar@{|->}[r]^(0.45){d'} \ar@{|->}[u]^{d"} & d'(y) & }
\]
In order to get the $E_2$ tableau, we have to take cohomology with respect to $d"$ and then $d'$. Let us compute $d_2^{0,1}([x])$ for $[x] \in E_2^{0,1}$ represented by an element $x \in FJ^{0,1}$, which therefore must satisfy $d"(x) = 0$ and $d'(x)$ is a $d"$ boundary, say $d'(x) = d"(y)$ with $y \in FJ^{1,0}$. Then 
\[d(x) = d'(x) = d"(y) = d(y) - d'(y) \]
which is cohomologous to $-d'(y)$ in the subquotient $E_2^{2,0}$ of $FJ^{2,0}$. So finally $d_2^{0,1}([x]) = [-d'(y)]$. 
\end{az}

\begin{az} Let us implement the above computation in our case. The complex
\[ k^\alg(X)^\ast \xrightarrow{\divisor} \Div(X \otimes k^\alg) \]
of the divisor map computes the relevant part of $\R\Gamma(X\otimes k^\alg,\Gm)$ and serves as $GI^\bullet$ with differential $d_{GI}$. In order to get an Eilenberg--MacLane resolution we take standard continuous cochains for $\Gal(k^\alg/k)$, the differential of which we write as $d^\vee$. The sign change $d" = (-1)^{p}d_{GI}$ guarantees that the two differentials of the double complex anticommute.
 
A Galois invariant line bundle $L \in \uPic_X(k)$ is represented by a divisor $D \in \Div(X \otimes k^\alg)$, that is linearly equivalent to its Galois conjugates, so there are functions $f_\sigma$ with 
\[\divisor(f_\sigma) = \sigma(D) - D\]
for each $\sigma \in \Gal(k^\alg/k)$.  The differential $d'(D) = d^\vee(D)$ is the cocycle $\sigma(D) - D$ which equals $\divisor(f_\sigma) = -d"(f_\sigma)$. We conclude, that 
\[d_2^{0,1}(L) = d^\vee(f_\sigma) = \sigma(f_\tau) f_{\sigma\tau}^{-1} f_\sigma. \]
In order to compute the Brauer obstruction $b(L)$ we need an actual line bundle representing $L$. We choose $\cL = \OO(D)$. The isomorphisms $\ph_\sigma : ^\sigma \cL \to \cL$ are identified with sections 
\[ \ph_\sigma \in \rH^0\big(X \otimes k^\alg, \cL \otimes (^\sigma\cL)^{-1}\big).\]
So $\divisor(\ph_\sigma) + D - \sigma(D) \geq 0$ and $\ph_\sigma = f_\sigma$ is a possible choice, hence
\[  b(L) = d^\vee(\ph_\sigma) = d^\vee(f_\sigma) = d_2^{0,1}(L),\]
as desired.
\end{az}

\end{appendix}

\bigskip

\noindent \textbf{Acknowledgments}

\noindent The paper was written while the author enjoyed visiting the Department of Mathematics of the University of Pennsylvania. I am grateful to Florian Pop for stimulus and discussions that lead to this paper, and his interest in my work in general. Thanks also go to the members of the Galois seminar at UPenn, especially Asher Auel and Ted Chinburg.








\end{document}